\documentclass[a4paper, 11pt]{amsart}
\pdfoutput=1
\usepackage{amssymb,amsmath,txfonts,mathrsfs,titletoc,appendix}
\usepackage{graphicx}
\usepackage{latexsym}
\usepackage[alphabetic]{amsrefs}
\usepackage{geometry}
\usepackage{fancyhdr}

 \newtheorem{thm}{Theorem}[section]
 \newtheorem{cor}[thm]{Corollary}
 \newtheorem{lem}[thm]{Lemma}

 \newtheorem{defn}[thm]{Definition}
 \newtheorem{rem}[thm]{Remark}
 \numberwithin{equation}{section}
\newtheorem{lem*}{Lemma}
\newtheorem{cor*}{Corollary}

\newenvironment{prooff}{\medskip \noindent
{\bf Proof.}}{\hfill \rule{.5em}{1em}
\\}

\begin{document}

\title{On extending calibrations}

\author{Yongsheng Zhang}
\address{{\it Current address:} School of Mathematics and Statistics, Northeast Normal University, 5268 RenMin Street, NanGuan District, ChangChun, JiLin 130024, P.R. China }
\email{yongsheng.chang@gmail.com}
\date{\today}

\begin{abstract}
This is a companion note of \cite{Z11} where the extension of local calibration pairs of smooth submanifolds is discussed.
Here we emphasize on the case of singular submanifolds. More precisely, we study when a calibration pair around the singular set of a submanifold can extend
to a local calibration pair about the entire submanifold.
Based upon \cite{Z2} several examples of particular interests under  the view of calibrated geometry are considered.

\end{abstract}
\keywords{metric, homologically mass-minimizing current, calibration with singularities, extension of calibration pairs} \subjclass[2010]{Primary~53C38, Secondary~49Q15, ~28A75}

\maketitle
\titlecontents{section}[0em]{}{\hspace{.5em}}{}{\titlerule*[1pc]{.}\contentspage}
\titlecontents{subsection}[1.5em]{}{\hspace{.5em}}{}{\titlerule*[1pc]{.}\contentspage}

\section{Introduction}\label{Section1} 
It was known since \cite{BDG} that there exist area-minimizing hypercones, i.e., Simons cones, in certain Euclidean spaces.
Very shortly afterward, Lawson \cite{BL} and others added many other area-minimizing hypercones.
However it remained unknown for decades whether these hypercones can be realized as tangent cones at some singular points
of homologically area-minimizing singular {\it compact} hypersurfaces in Riemannian manifolds.

In 1999, N. Smale \cite{NS} constructed first such examples through techniques of geometric analysis.
In this paper we shall show how to obtain such creatures thru the theory of calibrations.
We recall some definitions and the fundamental theorem of calibrated geometry (FTCG) in \S\ref{P}.
Then 
we establish our main result Theorem \ref{conecal}.
As applications examples are constructed in \S\ref{exs}.

Example 1 tells us how to make use of Theorem \ref{conecal} to create a metric such that
a homologically nontrivial compact singular hypersurface $(S,\mathscr S)$ (defined in \S\ref{P}) 
with ``nice" local behavior around $\mathscr S$ becomes homologically area-minimizing.
In our construction, the hypersurface  will be calibrated by a coflat calibration with exactly the same singular set $\mathscr S$. 

Example 2a shows that, based on a result of \cite{Z2}, one can find many 
examples of homologically nontrivial compact singular hypersurfaces
with the required behavior around singularities.
Combined with the construction of Example 1, quite a few examples similar to Smale's can be created by the theory of calibrations.

Example 2b conveys that there can exist a homologically area-minimizing smooth hypersurface
in some Riemannian manifold which however cannot be calibrated by any smooth calibration.
This disproves in some sense the reverse of the FTCG.

Example 2c discusses how to apply the method in Example 1 for homologically area-minimizing submanifolds of higher codimension.
Allowed local models of an isolated singular point consist of, for instance, all special {Lagrangian} cones.
(See
Joyce \cite{Joyce}, McIntosh \cite{Mc}, Carberry and McIntosh \cite{CM}, 
Haskins \cite{Haskins}, 
Haskins and Kapouleas \cite{HK}, \cite{HK2}, \cite{HK3} and etc. for the diversity.)

The case of non-orientable hypersurfaces is studied as well.
Example 3 relates to Murdoch's theory of twisted calibrations \cite{TM}.
In Example 4 we create a non-orientable compact singular hypersurface 
which is mass-minimizing in its homology class of integral currents mod 2 (introduced by Ziemer \cite{Ziemer}).

{\ }\\{\ }\\
$\text{\sc Acknowledgement}.$
This paper is an expansion of part of the author's Ph.D. thesis at Stony Brook University.
He is deeply indebted to his advisor Professor H. Blaine Lawson, Jr. for encouragement and guidance.
He also would like to thank Professor Frank Morgan for invaluable suggestions and Professor Brian White for helpful comments.
Part of the work was polished during the author's visit to the MSRI in Fall 2013.
{\ }\\{\ }

\section{Preliminaries}\label{P} 


\subsection{Calibrated Geometry}\label{cg}

We briefly recall definitions and notions as well as the FTCG. For details readers are referred to \cite{HL2}.

        \begin{defn}\label{calibration}
        A smooth form $\phi$ on a Riemannian manifold $(X,g)$ is called a $\mathrm{\mathbf{calibration}}$ if 
        $\sup_{X}\|\phi\|_{g}^*= 1$
        and
        $d\phi=0.$
        The triple $(X,\phi,g)$ is called a $\bold{calibrated}$ $\bold{manifold}$. 
        \end{defn}
        
        \begin{defn}
        
        By a $\mathrm{\mathbf{singular\ submanifold}}$ $(S,\mathscr S)$ with singular set $\mathscr S$, we mean a pair of closed
subsets $\mathscr S\subset S$ of X, such that $S-\mathscr S$ is an $m$-dimensional submanifold
and the Hausdorff $m$-measure $\mathcal H^m(\mathscr S) = 0$.
        \end{defn}

        \begin{defn}
        If $(S,\mathscr S)$ is an oriented (singluar) submanifold
        with $\phi|_{S-\mathscr S}$ equals to the volume form of $S-\mathscr S$,
        then
        $(\phi,g)$ is a $\mathrm{\mathbf{calibration\ pair}}$ of $S$ on $X$.
        We say $\phi$
        $\mathrm{\mathbf{calibrates}}$ $S$ and $S$ $\mathrm{\mathbf{can\ be\ calibrated}}$ in $(X,g)$.
        \end{defn}

When an $m$-dimensional current $T$ has local finite mass, we have decomposition
$$T = \overrightarrow T\cdot \|T\|\ a.e.\ \|T\|$$
where 
the {\it Radon} measure $\|T\|$ is characterized by
$ \int_X f\cdot d\|T\|=\sup\{T(\psi): \|\psi\|_{x,g}^*\leq f(x) \}$
for any nonnegative continuous function $f$ with compact support on $X$,
and $\overrightarrow T$ is a $\|T\|$ measurable tangent $m$-vector field a.e. with 
vectors $\overrightarrow T_x \in \Lambda^m T_xX$ of unit length in the dual norm of the comass norm.

\begin{defn}\label{calibratable}
Let $\phi$ be a calibration on $(X,g)$.
Then a current $T$ of local finite mass is $\bold{calibrated}$ by $\phi$, if 
$\phi_x(\overrightarrow T_x)=1\ a.a.\ x\in X\ \text{for}\ \|T\|.$
\end{defn}

\begin{rem}\label{coneccc}
Assume $S$ is a submanifold with only one singular point $p$ and $C_p$ is a tangent cone of $S$ at $p$.
Then the current $[[S]]=\int_S\cdot\ $ is calibrated by $\phi$ if and only if $S$ is calibrated, and moreover
either implies 
$\phi_p$ calibrates $C_p$ in $(T_pX,g_p)$.
\end{rem}

The following is the fundamental theorem of calibrated geometry.

\begin{thm}[Harvey and Lawson]\label{hl}
If $T$ is a calibrated current with compact support in $(X,\phi,g)$ and $T'$ is any compactly supported current homologous to $T$(i.e., $T-T'$ is a boundary and in particular $dT=dT'$), then
\begin{equation*}
\mathrm{\mathbf{M}}(T)\leq  \mathrm{\mathbf{M}}(T')
\end{equation*}
with equality if and only if $T'$ is calibrated as well.
\end{thm}

\begin{defn}
Let $\mathbf{spt}(f)$ be the support of $f$ where $f$ is a function.
For a current $T$, let $U_T$ stand for the largest open set with $\|T\|(U_T)=0$.
Then the support of $T$ is denoted by $\mathbf{spt}(T)=U_T^c$.
\end{defn}
{\ }

\section{Extension Results}\label{CWS}

\begin{thm}\label{conecal}
Suppose $S$ is an $m$-dimensional oriented compact submanifold with only one singular point $o$ in $(X^n,g)$
and it represents a nonzero class in the $\mathbb{R}$-homology of $X$.
 If a local part
 $B_{\epsilon}(o;g)\cap S$ for some $\epsilon>0$
 can be calibrated in the $\epsilon$-ball $(B_\epsilon(o;g),g)$, then
 there exists a metric $\hat g$
 coinciding with $g$ on $B_{{\frac{\epsilon}{2}}}(o;g)$ such that
 $S$ can be calibrated by a smooth calibration in $(X,\hat g)$.
\end{thm}

 \begin{figure}[ht]
\begin{center}
\includegraphics[scale=0.25]{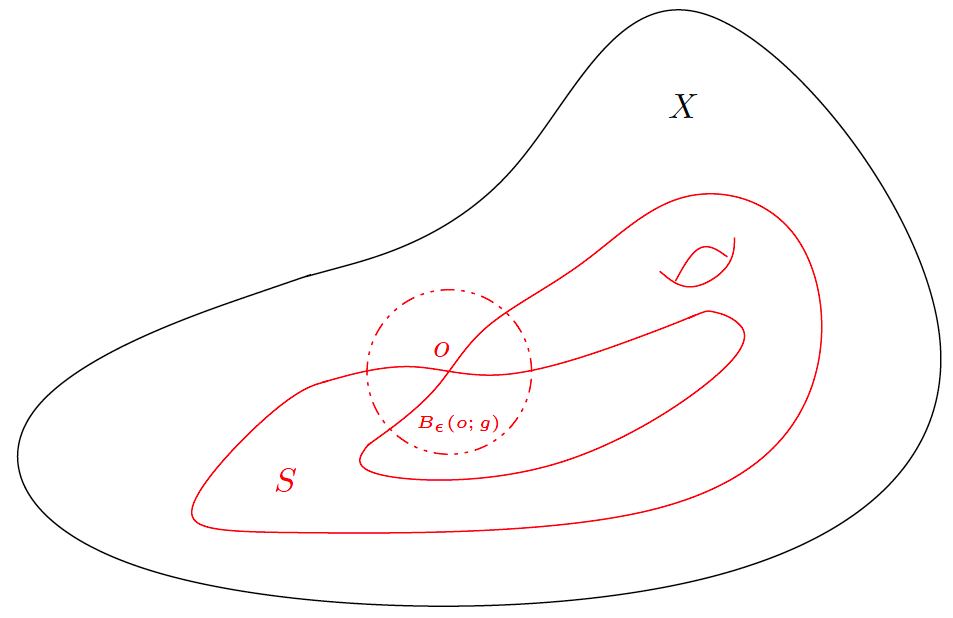}\ \ \ \ \ \ \ 
\end{center}
\end{figure}

\begin{rem}
In the theorem, 
$\frac{\epsilon}{2}$ can be replaced by $\varkappa\epsilon$ for any $0<\varkappa<1$.
\end{rem}

\begin{prooff}
 Suppose $\epsilon$ in the assumption is sufficiently small so that the open disc $D=B_\epsilon(o;g)$ corresponds to an open disc in some local chart. 
 Let $\tilde\phi$ be the local calibration form.
Then $\tilde\phi=d\psi$ where $\psi$ is some smooth $(m-1)$-form defined on $D$. Suppose the compact region $\Gamma_1\bigcup \Omega\bigcup\Gamma_2$
(given in the picture by the fiber structure induced by $g$ over the set $(\Gamma_1\bigcup \Omega\bigcup\Gamma_2)\cap S$ for small $h$)
is contained in $D-B_{{\frac{2\epsilon}{3}}}(o;g)$.
 \begin{figure}[ht]
\begin{center}
\includegraphics[scale=0.25]{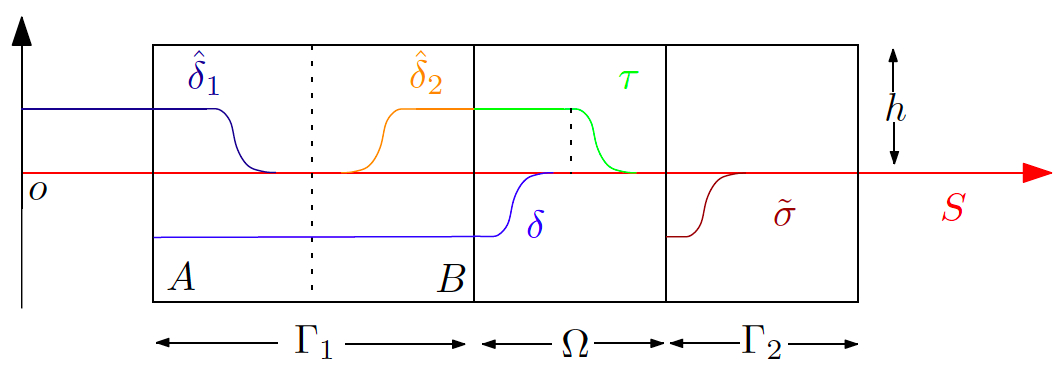}
\end{center}
\end{figure}

 Then $\pi_g^*\omega= d(\pi_g^*(\psi|_S))$ in $\Gamma_1\bigcup \Omega\bigcup\Gamma_2$
 where $\omega$ is the volume form of $S\bigcap(\Gamma_1\bigcup \Omega\bigcup\Gamma_2)$. 
 Set
 \[
 \Phi\triangleq d(\tau\psi+(1-\tau)\pi_g^*(\psi|_S))
 \]
 where $\tau$ is a cut-off function in $\Omega$ showed in the picture
 with value one near $\Gamma_1$ and zero near $\Gamma_2$. (The picture here is just an illustration, since the region ``hight" $h$ is generally smaller than one to guarantee
 no overlapping.)
 Since $\Phi(\overrightarrow{T_xS_g})
 =1$
 where $x\in S\bigcap(\Gamma_1\bigcup \Omega\bigcup\Gamma_2)$ and $\overrightarrow{T_xS_g}$ is the unique oriented unit horizontal $m$-vector 
 through $x$ on $S$,
 it can be achieved by shrinking $h$ (with respect to $g$)
  that the smooth function
 $\Phi(\overrightarrow{T_yS_g})>\frac{1}{2}\text{ on }\Gamma_1\bigcup \Omega\bigcup\Gamma_2$
 where $y$ in $\Gamma_1\bigcup \Omega\bigcup\Gamma_2$
 and $\overrightarrow{T_yS_g}$ is the unique oriented unit horizontal
 (to the disc-fibration $\mathscr F$ generated by the exponential map restricted to normal directions along
 $S\bigcap(\Gamma_1\bigcup \Omega\bigcup\Gamma_2)$)
 $m$-vector at $y$ with respect to $g$.
 Set 
 \[
 \bar g=f\cdot g
\text{ where } 
f=\delta+(1-\delta)(\Phi(\overrightarrow{T_yS_g}))^{\frac{2}{m}}
\]
on $\Gamma_1\bigcup \Omega\bigcup\Gamma_2$.
$\Phi=\tilde\phi$ on $\bold{spt}(\delta)$. Since $(\tilde\phi,g)$ is a local calibration pair given in the assumption,
 we know $f\geq (\Phi(\overrightarrow{T_yS_g}))^{\frac{2}{m}}$ on $\Gamma_1\bigcup\Omega\bigcup \Gamma_2$
 and $f\equiv1$ on $\Gamma_1$.
 Then $\bar g$ can extend on $\Upsilon$, the region embraced by the ``curve" in the picture below (which is an ``$h$-disc bundle'' containing $\Gamma_1\bigcup\Omega\bigcup\Gamma_2$),
such that
\\
(a). $\Phi$ (naturally extended on $\Upsilon$) calibrates $S\bigcap(\Upsilon-\Omega)$ in $(\Upsilon-\Omega,\bar g)$,
\\
(b). $\bar g=g$ in $\Gamma_1$, and
\\
(c). $\Phi(\overrightarrow{T_yS_{\bar g}})\leq 1$ on $\Upsilon$ with equality on $\Upsilon-\Gamma_1-\Omega$,
where $\overrightarrow{T_yS_{\bar g}}$ is the unique oriented unit horizontal (to $\mathscr F$) $m$-vector at $y$ with respect to $\bar g$.

 \begin{figure}[ht]
\begin{center}
\includegraphics[scale=0.25]{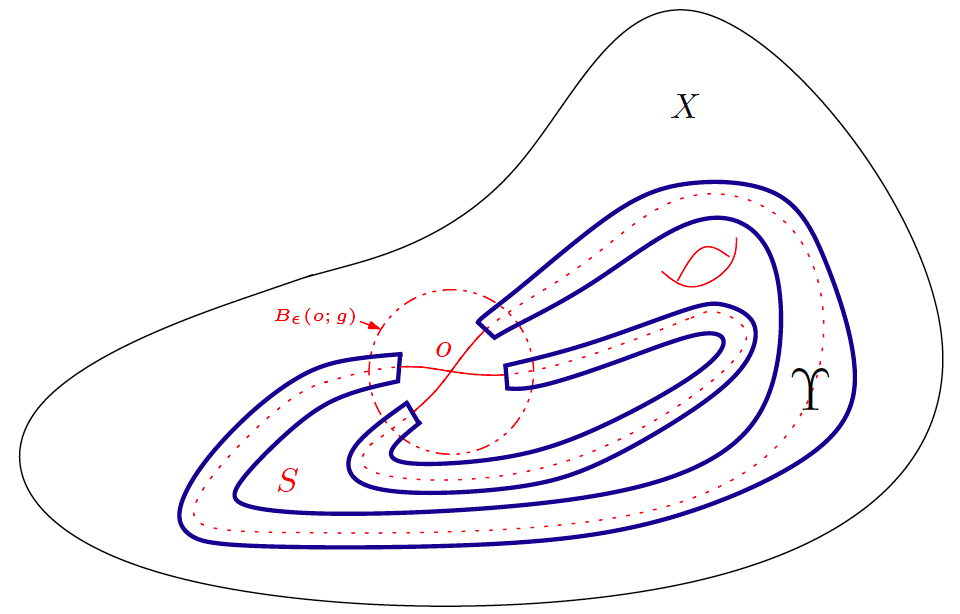}
\end{center}
\end{figure}

 In order to glue $\bar g$ and $g$ together and meanwhile to guarantee $\Phi$ a calibration,
 we need the following
 powerful lemmas from \cite{HL1}. 
 \begin{lem}[Harvey and Lawson] \label{hl1}
 Let $\xi\in \Lambda^p \mathbb{R}^n$ be a simple p-vector with
  $V = span\{\xi\}$. Suppose $\phi\in\Lambda^p\mathbb{R}^n$ satisfies  $\phi(\xi)= 1$. 
  Then there exists a unique oriented complementary subspace $W$ to $V$ 
  with the following property.
   For any basis $v_1, \cdots, v_n$ of $\mathbb{R}^n$ such that $\xi=v_1\wedge... \wedge v_p$ 
   and $v-{p+1}, \cdots, v_n$ is basis for $W$, one has that
  \begin{equation}
 \phi=v_1^*\wedge \cdots \wedge v_p^*+ \sum a_Iv_I^*,
 \end{equation}
 where $a_I=0$ whenever $i_{p-1}\leq p$. 
 \end{lem}
 \begin{lem}
 [Harvey and Lawson]
 \label{hl2}
Let $\phi,\  V=	span \{\xi\}$, and $W$ 
be as in Lemma \ref{hl1}.
Consider an inner product $<\cdot, \cdot>$
on $\mathbb{R}^n$ such that $V\perp W$ and $\|\xi\|=1$.
Choose any constant $C^2 > 
\bigl( \begin{smallmatrix} n\\ p\end{smallmatrix} \bigr)
\|\phi\|^*$
and define a new inner product on $\mathbb{R}^n=V\oplus W$
by setting $<\cdot,\cdot>'=<\cdot,\cdot>_V+C^2<\cdot,\cdot>_W$.
Then in this new metric we have
$$\|\phi\|^*=1\ and\ \phi(\xi)=\|\xi\|=1.$$
 \end{lem}
 \begin{rem}\label{imp}
If 
 $\phi(\xi)=\vartheta$ (positive) instead of one,
 one can apply {\it Lemma \ref{hl1}} to $\vartheta^{-1}\phi$ to get a similar conclusion that
 $\|\phi\|^*=\vartheta,\ \|\xi\|=1\ and\ \phi(\xi)=\vartheta$
 by choosing $C^2>\vartheta^{-1} \bigl( \begin{smallmatrix} n\\ p\end{smallmatrix} \bigr)
\|\phi\|^*$.
 \end{rem}
 By applying {Lemma \ref{hl1}} to $\Phi$, $\overrightarrow{T_yS_{\bar g}}$ and
 $\bar g$ on $\Upsilon$, one can get a smoothly varying $(n-m)$-dimensional 
 plane field $\mathscr W$ transverse to the horizontal directions in $\Upsilon$.
 Following {Lemma \ref{hl2}}, {Remark \ref{imp}} and property (c),
 for any metric $g_{\mathscr W}$ along $\mathscr W$,
 there exists a sufficiently large constant $\bar \alpha$ (due to the compactness of $\Upsilon$) such that,
 under $\tilde g= \bar g^h\oplus\bar \alpha g_{\mathscr W}$ on $\Upsilon$, 
\[
 \|\Phi\|_{\tilde g}^*=
\Phi(\overrightarrow{T_yS_{\bar g}})\leq 1.
\]

Now construct a smooth metric $\check g$ on $\Xi$ as follows based on property (b).
$$
\check g=
\begin{cases}\label{aaa}
g & \text{near}\ o\\
g+(1-\hat \delta_1)((0\cdot {\bar g^h})\oplus \bar \alpha g_{\mathscr W}) & \text{on} \ A\\
(1-\hat \delta_2)((0\cdot g^h)\oplus g^\nu) +\tilde g & \text{on} \ B\\
\tilde g & \text{on}\ \Omega\\
\tilde \sigma
\tilde g+(1-\tilde \sigma)
\bar g & \text{on}\ \Gamma_2\\
\bar g & \text{far away from}\ o\\
\end{cases}
$$
Here $\oplus$ means the orthogonal splitting of a (pseudo-)metric and $+$ is the usual addition between two (pseudo-)metrics.

On $\Gamma_2$,  ${\mathscr W}$ is exactly the distribution of fiber directions of $\mathscr F$
and $\Phi=\pi_g^*(\omega)$ is a simple horizontal $m$-form.
By {Lemmas 2.13}, {2.14} and {2.15} of \cite{Z11}, $\Phi$ is a calibration in $(\Xi,\check g)$.
Note that $\Xi$ can be retracted to $S$ trough a strong deformation retraction.
Therefore applying the gluing tricks of forms and metrics in \cite{Z11} on smaller regions of $S$
produces a global calibration pair $(\hat\Phi, \hat g)$ of $S$.
 \begin{figure}[ht]
\begin{center}
\includegraphics[scale=0.25]{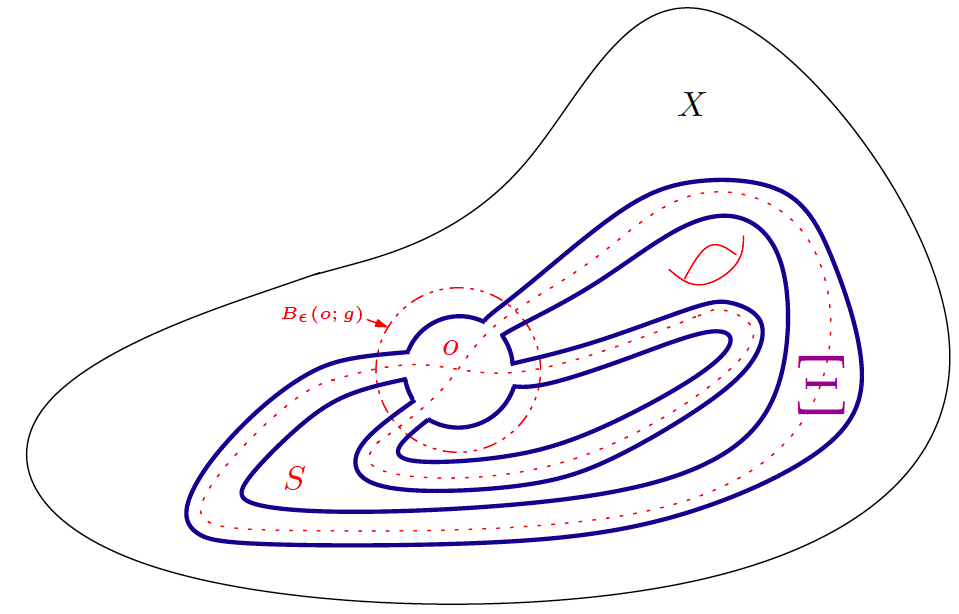}
\end{center}
\end{figure}
\end{prooff}

            By noticing that the comass function of a smooth form of codimension one is smooth we get the following refinement.
            \begin{cor}\label{codim1}
Suppose $S$ is an oriented compact hypersurface with only one singular point $o$ in $(X,g)$
and it represents a nonzero class in the $\mathbb{R}$-homology of $X$.
 If $B_{\epsilon}(o;g)\cap S$ for some $\epsilon>0$ can be calibrated by a calibration singular only at $o$ in the $\epsilon$-ball $(B_\epsilon(o;g),g)$, then
 there exists a metric $\hat g$ in the conformal class of $g$
 coinciding with $g$ on $B_{{\frac{\epsilon}{2}}}(o;g)$ such that
 $S$ can be calibrated by a calibration singular only at $o$ in $(X,\hat g)$.
\end{cor}

In fact it does not have to require that $S$ is a retract of some open neighborhood of $S$ for the last step of the proof. 
Whenever there exists a global defined form
which represents $[\Phi]$ on some open neighborhood of $S$, our construction applies. 
\begin{thm}
Suppose $(S,\mathscr S)$ is an $m$-dimensional oriented connected compact singular submanifold in $(X,g)$
with $[S]\neq[0]$ in $H_m(X;\mathbb R)$.
Assume $V\bigcap S$
where $V$ is an open neighborhood of $\mathscr S$
can be calibrated in $(V,g|_V)$.
If 
\[
i^*: H^m(X;\mathbb R)\rightarrow H^m(B_{\epsilon}(S;g);\mathbb R),
\]
where $B_{\epsilon}(S;g)$ is the $\epsilon$-neighborhood of $S$ under $g$,
is onto for a sufficiently small positive $\epsilon$,
 then there exists a metric $\hat g$
 coinciding with $g$ in $B_{{\frac{\epsilon}{2}}}(\mathscr S;g)$ such that
 $S$ can be calibrated in $(X,\hat g)$.
\end{thm}
               
                \begin{rem}
                By {\it Almgren}'s big regularity theorem, $\mathscr S$ has codimension at least $2$ in $S$.
                By $\bold{spt}(d[[S]])\subseteq\mathscr S$, $d[[S]]=0$ and therefore $[S]$ makes sense.
                \end{rem}

                \begin{rem}
When $\mathscr S$ is a smooth submanifold, for a sufficiently small $\epsilon>0$, $B_{\epsilon}(S;g)$ can be strongly retracted to $S$.
\end{rem}
{\ }

\section{Further Applications}\label{exs}
Sometimes it is useful to consider calibrations with singularities.
By \cite{Z2} every area-minimizing cone $C_{n,m}\subset\mathbb R^{n+m+2}$ enjoys a calibration singular only at the origin.
\\

\textbf{Example 1:} 
When the local model around $o$ in Theorem \ref{conecal} is a
{\em Simons} cone over $S^{r-1}\times S^{r-1}$ with $r\geq 4$,
one has an $SO(r)\times SO(r)$ invariant smooth calibration $\tilde \phi$ on $\mathbb R^{2r}-\{0\}$.
Follow the proof of Theorem \ref{conecal} to get $\Phi$ on $\Xi-o$ and $\check g$ on $\Xi$. 
By {\em Mayer-Vietoris} sequence,
$$H^{2r-2}(S^{2r-1})\rightarrow H^{2r-1}(\Xi)\rightarrow H^{2r-1}(\Xi-o)\rightarrow H^{2r-1}(S^{2r-1})\rightarrow H^{2r}(\Xi),$$
where $S^{2r-1}$ is a small sphere centered at $o$.
Since
\[
\int_{S^{2r-1}}\Phi=0
\]
and
 there exists a strong deformation retraction from $\Xi$ to $S$,
 one can obtain a smooth form $\check \phi$ on $X$ such that
 \[
 \check \phi|_{\Xi-o}-\Phi=d\check \psi
 \]
 for some smooth $(2r-2)$-form $\check \psi$ on $\Xi-o$.
Then, away from $S$, glue $\check \phi$ and $\Phi$ together to get a smooth form $\hat \Phi$ on $X-o$,
and meanwhile extend $\check g$ to $\hat g$
making $\hat\Phi$ a calibration on $X-o$
(c.f. \cite{Z11}).

Due to Theorem \ref{hl},
$[[S]]$ is mass-minimizing in its current homology class.
However, it is impossible to calibrate $S$ by any smooth calibration $\bar \Phi$ on $(X,\hat g)$ (actually for any metric). Since if so, according to Remark \ref{coneccc}
the tangent cone of $S$ at $o$, a {\em Simons} cone, would be calibrated in $(T_oX, \bar\Phi_o,\hat g_o)$. But $\bar\Phi_o$ can calibrate certain hyperplanes only. Contradiction!
\\{\ }

\textbf{Example 2a:} On an oriented compact $(2r-1)$-dimensional smooth manifold $T$,
take a small disc neighborhood and trivially embed $K=S^{r-1}\times D^r$ into it
(one can of course use other embedding of $S^{r-1}\hookrightarrow T$  as long as the normal bundle is trivial).
After surgery along $S^{r-1}\times S^{r-1}$, denote the obtained manifold by $T'$.
Then $T$ and $T'$ are cobordant through some oriented $2r$-dimensional smooth manifold $W$.
Specifically, $W$ can be taken as the union of $[-0.5,0.5]\times (T-K)$ and the region between $\{f=-0.5\}$ and $\{f=0.5\}$ in the picture
via a family of gluing diffeomorphisms of $S^{r-1}\times S^{r-1}$. ($K$ is identified with $\{f=-0.5\}$.)
Here $f$ is defined on the unit ball in $\mathbb R^r\times \mathbb R^r$
given by 
$f(\overrightarrow x, \overrightarrow y)=-\| \overrightarrow x\|^2+\|\overrightarrow y\|^2.$
 \begin{figure}[ht]
\begin{center}
\includegraphics[scale=0.25]{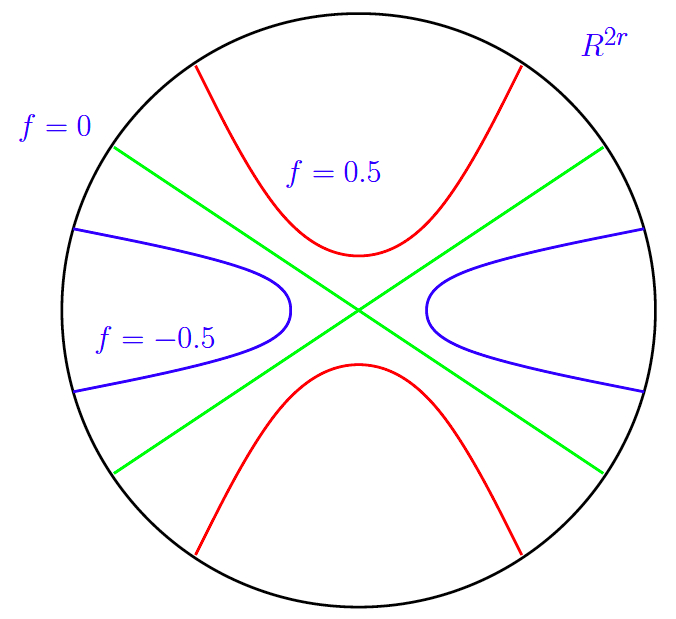}
\end{center}
\end{figure}
Observe that $f^{-1}(0)$ is a truncated {\em Simons} cone and the foliation is similar to the situation around a saddle point.

Take two copies of $W$. Glue the same boundaries and one gets an orientable compact $2r$-dimensional manifold $X$.
Now extend the Euclidean metric on the region between $\{f=-0.25\}$ and $\{f=0.25\}$ in the first copy to a metric on $X$.
Let $S$ be the slice corresponding to $f=0$.
Apparently $[S]\neq[0]$ in $H_{2r-1}(X;\mathbb R)$ (by intersection number method).
Then the above arguments show that 
$S$ can be calibrated by a calibration $\Phi$ singular only at the origin with respect to some metric $g$ on $X$.
\begin{rem}\label{88}
By cross-products examples with more complicated singularity can be generated.
For instance, 
$S\times S$ with singularity $S\vee S$ is calibrated by a coflat calibration
with singular set 
$
S\vee S
$ 
in the cartesian product $(X, g)\times (X, g)$.
\end{rem}
{\ }

\textbf{Example 2b:} 
In Example 2a, choose a smooth ``fiber", for example,
$M\triangleq(T-K)\bigcup \{f=-0.3\}$
in the first copy of $W$. Note that $\Phi$ is already a coflat calibration of $S$ on $(X, g)$. According to the method in \cite{Z11}, one can modify the calibration to $\tilde \Phi$ and conformally change $g$ to $\tilde g$ in a neighborhood of $M$ away from $S$
such that $\tilde \Phi$ becomes a coflat calibration calibrating both $S$ and $M$ in $(X, \tilde g)$. 

However the homologically mass-minimizing smooth submanifold $M$ cannot be calibrated by any smooth calibration in $(X, \tilde g)$. If it were,  then $S$ must be calibrated by the same smooth calibration as well which would lead to a contradiction as before.
This implies that all the coflat calibrations of $M$
in $(X,\tilde g)$ share at least a common singular point.
For such creatures of higher codimension,
one can consider 
$M\times \{ \text{a point}\}$
in the Riemannian product of $(X,\tilde g)$ and a compact oriented manifold.
\\

\textbf{Example 2c:}
Generally, the construction in Example 2a applies to hypercones only.
For cones of higher codimension
we introduce the following construction.
Suppose $C\subset \mathbb R^n$ is a $k$-dimensional cone which can be calibrated by some calibration $\phi=d\psi$
possibly singular at the origin.
Then similarly as in \cite{NS}
consider $\Sigma_C\triangleq (C\times \mathbb R)\bigcap S^{n}(1)$ in $ \mathbb R^{n+1}$.
Choose an $n$-dimensional oriented compact manifold $T$
with nontrivial $H_k(T;\mathbb R)$.
Let $M$ be an embedded oriented connected compact submanifold with $[M]\neq [0]\in H_k(T;\mathbb R)$.
Taking smooth disks around a point of $M$ and a smooth point of $\Sigma_C$ respectively
one can simultaneously connect $T$ and $S^n(1)$, $M$ and $\Sigma_C$
through one surgery along $S^0\times S^n$ (i.e., connect sum).
Denote by $X$ and $S$
the resulted manifold and submanifold (singular at two points).
Then $[S]\neq 0\in H_k(X;\mathbb R)$
and there exists a calibration pair of $S$ on $X$
according to Theorem \ref{conecal}.
\\{\ }

Next we consider the non-orientable case.
\\

\textbf{Example 3:}
One can do blowing-ups of $X$ at several points away from the singular point $o$ of $S$ in Example 2a.
An interesting simple case is to blow up at a smooth point of $S$.
Call the resulted manifold and submanifold $\check X$ and $\check S$.
Note that $\check S$ inherits a natural orientation from $M$.
Moreover, $[\check S]\neq 0$.  This can be seen either by pairing with a suitable closed form
or by lifting it to the double cover $\overline X$ of $\check X$ corresponding to $\pi_1\mathbb RP^{2r}$.
By the first idea, one can get some calibration pair as before.
The second can generate
some $\mathbb Z_2$-invariant calibration pair on the double cover, which induces a calibration pair on $\check X$.
Then under the resulted metric
$\check S$ can be calibrated by some coflat calibration singular only at $o$ in the non-orientable $\check X$.

\begin{figure}[htbp]
\begin{minipage}[c]{0.3\textwidth}
  \includegraphics[scale=0.11]{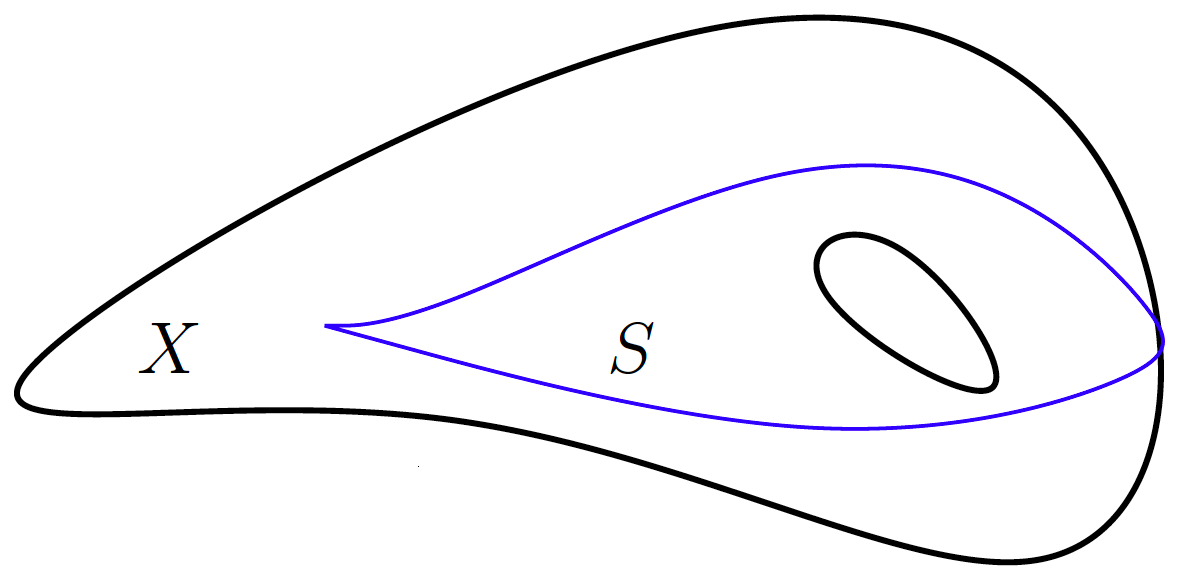}
\end{minipage}%
\ \ \ \ \ \ \ \ \ \ \Huge $\leadsto$
\begin{minipage}[c]{0.4\textwidth}
  \includegraphics[scale=0.21]{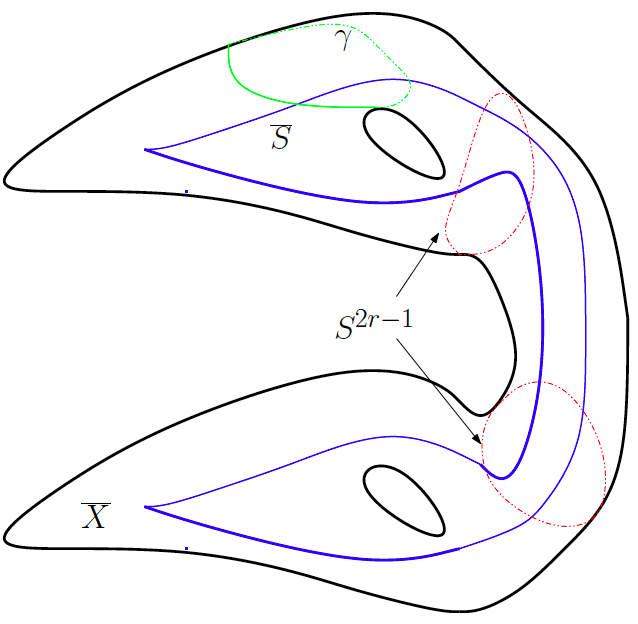}
\end{minipage}
\end{figure}
{\ }

\textbf{Example 4:}
Based upon $C_{3,4}$ one can get
an eight-dimensional oriented compact connected submanifold
$S$ with one
singular point in 
some oriented manifold $X^9$
with $[S]\neq [0]\in H_8(X;\mathbb R)$
by the method of Example 2a.
Now blow up at a regular point of $S$.
Call the resulted manifold and submanifold $\check X$ and $\check S$ respectively.

By the {\em Seifert-van\ Kampen} theorem $\pi_1(\check X)\cong\pi_1(X)*\pi_1(\mathbb RP^8).$
Similarly, the isomorphism of $\pi_1(\mathbb RP^8)\cong \mathbb Z_2$
trivially extends to a homomorphism $:$
$\pi_1(\check X)\rightarrow \mathbb Z_2$,
which canonically determines a two-sheeted cover $\overline X$ of $\check X$.
Denote the lifting of $\check S$ by $\overline S$.
Note that $\overline X\cong X\# X$ and $\overline S\cong S\#S^{\text{opposite orientation}}$.
By {\em Mayer-Vietoris} sequences
$$H_8(\overline X;\mathbb Z)\cong H_8(X;\mathbb Z)\oplus H_8(X;\mathbb Z), \text{and}$$
$$[\overline S]= [(S,-S)]\neq[0]\ \text{in} \ H_8(\overline X;\mathbb Z).$$
Now create a $\mathbb Z_2$-invariant metric $\bar g$ on $\overline X$ such that
the orientable $\overline S$ can be calibrated (by a twisted calibration in the sense of \cite{TM}).

Given a triangulation $\check S$ can be viewed as an integral current by assigning chambers local orientations.
Also note that $\check S$ induces a $d$-closed integral current mod $2$,
$[[\check S]]_2$ (see \cite{Ziemer}), representing a non-zero $\mathbb Z_2$-homology class $[\check S]_2$.
We want to show that $[[\check S]]_2$ is $\mathrm{\mathbf M}^2$-minimizing in $[\check S]_2$ under the induced metric
$\check g$ on $\check X$,
where $\bold{the\ mass}$ $\mathrm{\mathbf{M}}^2(\cdot)$ of an integral current mod $2$
is the infimum of the mass of all integral representatives.

Suppose $K-[[\check S]]_2=dW$ in the sense of mod 2 for an integral current $K$ of finite mass and $W$ a top dimensional integral current mod 2.
Then the lifting to $\overline X$ becomes $\overline K-[[\overline S]]=d\overline W$ in the sense of mod $2$.
(Since $\overline S$ is orientable, $[[\overline S]]$ is an integral current up to a choice of orientation.)
By $\overline X$'s being oriented and $\overline W$'s being of top dimension,
$\overline W$ comes from
the quotient of $\tilde W$ by $2$
where $\tilde W$ is the integral current with multiplicity one on $\mathrm{\bold{spt}}(\overline W)$
and orientation inherited from $\overline X$.
Restrict $\tilde W$ to the connected component of $\bold{spt}(\overline W)$ to $\overline S$
and
denote it by
$\tilde W^\circ$.
Assign $[[\overline S]]$ the orientation induced from $\tilde W^\circ$.  
Let $-\overline K^\circ\triangleq d\tilde W^\circ-[[\overline S]]$.
It follows
\[
\mathrm{\mathbf M}_{\check g}(\check S)=
\frac{1}{2}\mathrm{\mathbf{M}}_{\bar g}([[\overline S]])\leq\frac{1}{2}\mathrm{\mathbf{M}}_{\bar g}(\overline K^\circ)\leq \mathrm{\mathbf{M}}_{\check g}(K).
\]
Running $K$ through all the integral representatives of $[[\check S]]_2$ one has
\[
\mathrm{\mathbf M}_{\check g}(\check S)=\mathrm{\mathbf M}_{\check g}^2([[\check S]]_2).
\]
Let $K_2$ be the integral current mod $2$ of an integral current $K$ with $[K_2]=[\check S]_2$.
Then
\[
\mathrm{\mathbf M}_{\check g}^2([[\check S]]_2)\leq \mathrm{\mathbf M}^2_{\check g}(K_2).
\]
Namely, $[[\check S]]_2$ is $\mathrm{\mathbf M}_{\check g}^2$-minimizing (of mass $\mathrm{\mathbf M}_{\check g}(\check S)$) in its homology class.
{\ }
\\

\section{Questions}

We would like to end up this paper with two further questions.
\\
 
 \textbf{Question A:\ } Under the same hypotheses of Theorem \ref{conecal},
is it possible to {\it conformally} change $g$ such that
$S$ can be calibrated with respect to the new metric? (For the case of codimension larger than one.)
\\

\textbf{Question B:\ } 
Suppose a current $T$ is homologically mass-minimizing in a fixed ambient manifold $(X,g)$, is there some sufficient criterion for $T$ being calibrated by some smooth calibration?


{\ }
\\
{\ }
{\  }

\end{document}